\documentclass[12pt]{amsart}
\usepackage{amscd,amssymb}
\usepackage[arrow,matrix]{xy}
\usepackage[plainpages,backref,urlcolor=blue]{hyperref}

\theoremstyle{plain}
\newtheorem{thm}[subsection]{Theorem}

\theoremstyle{definition}
\newtheorem{rk}[subsection]{Remark}
\newtheorem{definition}[subsection]{Definition}
\newtheorem{ex}[subsection]{Example}

\numberwithin{equation}{section}
\newcommand{\OO}{{\mathcal O}}

\newcommand{\I}{{\mathcal I}}

\newcommand{\C}{\mathbb{C}}

\newcommand{\PP}{\mathbb{P}}

\newcommand{\N}{\mathbb{N}}

\DeclareMathOperator{\coker}{coker}

\DeclareMathOperator{\defect}{def}

\DeclareMathOperator{\indeg}{indeg}

\begin{document}

\title [Versality and bounds of global Tjurina numbers]
{Versality, bounds of global Tjurina numbers and logarithmic vector fields along hypersurfaces with isolated singularities}

\author[A. Dimca]{Alexandru Dimca$^1$}
\address{Universit\'e C\^ ote d'Azur, CNRS, LJAD and INRIA, France and Simion Stoilow Institute of Mathematics,
P.O. Box 1-764, RO-014700 Bucharest, Romania}
\email{dimca@unice.fr}

\thanks{$^1$ This work has been partially supported by the French government, through the $\rm UCA^{\rm JEDI}$ Investments in the Future project managed by the National Research Agency (ANR) with the reference number ANR-15-IDEX-01 and by the Romanian Ministry of Research and Innovation, CNCS - UEFISCDI, grant PN-III-P4-ID-PCE-2016-0030, within PNCDI III.}



\begin{abstract} We recall first the relations between the syzygies of the Jacobian ideal  of the defining equation for a projective hypersurface $V$ with isolated singularities and the versality properties of $V$, as studied by du Plessis and Wall. Then we show how the bounds on the global Tjurina number of $V$ obtained by du Plessis and Wall
lead to substantial improvements of our previous results on the stability of the reflexive sheaf $T\langle V\rangle$ of logarithmic vector fields along $V$, and on the
 Torelli property in the sense of Dolgachev-Kapranov of $V$. 

\end{abstract}

\maketitle


\section{Introduction} \label{sec:intro}

Let $V:f=0$ be a degree $d$ singular hypersurface in the complex projective space $\PP^n$, having only isolated singularities. Let $S=\C[x_0, \dots, x_n]$ be the graded polynomial ring, and consider the graded $S$-module of {\it Jacobian syzygies or Jacobian relations} of $f$ defined by
\begin{equation} 
\label{eqAR1}
AR(f)=\{(a_0, \dots,a_n) \in S^{n+1} \ : \ a_0f_0 + \dots + a_nf_n=0\},
\end{equation}
where $f_j$ denotes the partial derivative of the polynomial $f$ with respect to $x_j$ for $j=0, \dots, n.$
This module has a natural $S$-graded submodule $KR(f)$, the module of
{\it Koszul syzygies or Koszul relations} of $f$, defined as the submodule 
 spanned by obvious relations of the type $f_{j}f_{i}+(-f_{i})f_{j}=0$. Note that the syzygies in $AR(f)$ are regarded as vector fields annihilating $f$ in the papers by A. du Plessis and C.T.C. Wall, while the Koszul syzygies are called Hamiltonian vector fields.
 The quotient
\begin{equation} 
\label{eqER}
ER(f)=AR(f)/KR(f)
\end{equation}
is the graded module of {\it essential Jacobian relations}.  
If $K^*(f)$ denotes the (cohomological) Koszul complex of $f_{0},...,f_{n}$ see \cite{DBull}, then one clearly has
\begin{equation} 
\label{eqKR}
ER(f)=H^n(K^*(f))(-n).
\end{equation}
Using these two graded $S$-modules, we introduce two numerical invariants for the hypersurface $V:f=0$ as follows.
The integer
\begin{equation} 
\label{mdr}
mdr(f)=\indeg(AR(f))=\min\{k \ : \ AR(f)_k \ne 0 \}
\end{equation}
is called the {\it minimal degree of a relation} for $f$, while the integer
\begin{equation} 
\label{mder}
mder(f)=\indeg(ER(f))=\min\{k \ : \ ER(f)_k \ne 0 \}
\end{equation}
is called the {\it minimal degree of an essential relation} for $f$.
From the definition, it is clear that $mdr(f) \leq mder(f)$ with equality if
$mdr(f)<d-1$. Note also that $0 \leq mdr(f) \leq d-1$ and $0 \leq mder(f) \leq n(d-2)$, where the last inequality follows from \cite[Corollary 11]{CD}, see also Theorem \ref{thmA1} below.
It is also clear that $mdr(f)=0$ if and only if $V$ is a cone, case excluded in our discussion from now on.

Let $\alpha_V$ be the Arnold exponent of the hypersurface $V$, which is by definition the minimum of the Arnold exponents of the singular points of $V$, see \cite{DS2}. Using Hodge theory, one can prove that
\begin{equation} 
\label{van1}
   mder(f) \geq \alpha_Vd -n,
\end{equation}
under the {\it additional hypothesis that all the singularities of $V$ are weighted homogeneous}, see \cite{DS2}. 
This inequality is the best possible in general, as one can see by considering hypersurfaces with a lot of singularities, see \cite{DiSt2}. However, for situations where the hypersurface $V$ has a small number of singularities this result is far from optimal, and in such cases one has the following inequality
\begin{equation} 
\label{van2}
mder(f) > n(d-2)-\tau(V),
\end{equation}
where $\tau(V)$, the Tjurina number of $V$, is the sum of the Tjurina numbers of all the singularities of $V$, see \cite{DAlgGeo}.

Jacobian syzygies and these two invariants $mdr(f)$ and $mder(f)$ occur in a number of fundamental results due to A. du Plessis and C.T.C. Wall, see \cite{duP1,duPCTC99, duPCTC00,duPCTC01},
some of which we recall briefly below.
The first class of their results are devoted to the versality properties of projective hypersurfaces. These are recalled in section 2, where we explain that \cite[Theorem 1.1]{duP1}, which is stated as Theorem \ref{thmA} below, is essentially equivalent to the first part of \cite[Theorem 1]{DBull}, which is stated as Theorem \ref{thmA1} below
for the reader's convenience.

The second class of results by A. du Plessis and C.T.C. Wall are related to finding lower and upper bounds for the global Tjurina number $\tau(V)$. Their main result in this direction is \cite[Theorem 5.3]{duPCTC01}, which is stated as Theorem \ref{thmC} below.
We show that this result can be used to greatly strengthen two of our main results in \cite{DAlgGeo},  one on the stability of the reflexive sheaf $T\langle V\rangle$ of logarithmic vector fields along a surface $V$, and the other on the
 Torelli property in the sense of Dolgachev-Kapranov of  the hypersurface $V$, see Theorems \ref{thm0} and \ref{thm1} below. Since the proofs of our results given in \cite{DAlgGeo} are rather long and technical, we present here only the minor changes in these proofs, possible in view of du Plessis and Wall's result in Theorem \ref{thmC}, and leading to much stronger claims, as explained in Remarks \ref{rk0} and \ref{rk1}.

\section{Versality of hypersurfaces with isolated singularities} \label{sec2}

Fix an integer $a \geq 0$, and call the hypersurface $V:f=0$
{\it $a$-versal}, resp. {\it topologically $a$-versal}, if the singularities of $V$ can be simultaneously versally, resp. topologically versally,  deformed by deforming the equation $f$, in an affine chart $\ell \ne 0$ with $\ell \in S_1$, containing all of the singularities, by the addition of all monomials of degree $n(d-2)-1-a$. Otherwise, we say that $V$ is {\it (topologically) $a$-non-versal}. With the above notation, one has the following result proved by A. du Plessis, see \cite[Theorem 1.1]{duP1}.
\begin{thm}
\label{thmA}
The hypersurface $V:f=0$ is $a$-non-versal if and only if $$a \geq mder(f).$$
\end{thm}
Let $\Sigma$ be the singular subscheme of $V$, defined by the Jacobian ideal of $f$ given by
$$J_f=(f_0, \dots,f_n) \subset S.$$
 Then, for $p$ a singular point of $V$, one has an isomorphism of Artinian $\C$-algebras
$$\OO_{\Sigma,p}=T(g),$$
where $g=0$ is a local equation for the germ $(V,p)$ and $T(g)$ is the Tjurina algebra of $g$, which is also the base space of the miniversal deformation of the isolated singularity $(V,p)$. More precisely, one has
\begin{equation} 
\label{Tju}
T(g)=\frac{\C\{y_1,\dots,y_n\}}{J_g +(g)},
\end{equation}
where $(y_1,\dots,y_n)$ is a local coordinate system centered at $p$
and $J_g$ is the Jacobian ideal of $g$ in the local ring $\OO_{\PP^n,p}=\C\{y_1,\dots,y_n\}$. Note that, for any integer $k$, one can consider the natural evaluation morphism
$$eval_k:S_k \to \oplus _{p\in \Sigma} \OO_{\Sigma,p}, \  \  h \mapsto ([h/\ell^k]_p)_{p\in \Sigma},$$
computed in the chart $\ell \ne 0$. Alternatively, $eval_k$ is the morphism
$$H^0(\PP^n, \OO_{\PP^n}(k) )\to H^0(\PP^n, \OO_{\Sigma}(k))=H^0(\PP^n, \OO_{\Sigma}),$$
induced by the exact sequence
\begin{equation} 
\label{ES1}
0 \to \I_{\Sigma} \to \OO_{\PP^n} \to \OO_{\Sigma} \to 0,
\end{equation}
where $\I_{\Sigma}$ is the ideal sheaf defining the singular subscheme $\Sigma$.
We set
\begin{equation} 
\label{defdef}
\defect_k(\Sigma)=\dim (\coker eval_k)=\dim H^1(\PP^n, \I_{\Sigma}(k)),
\end{equation}
the {\it defect of $\Sigma$ with respect to homogeneous polynomials of degree} $k$.
It follows that
the hypersurface $V:f=0$ is $a$-versal if and only if the corresponding evaluation morphism
$eval_{n(d-2)-1-a}$ is surjective, i.e. the defect
$\defect_{n(d-2)-1-a}(\Sigma)$ vanishes. We see in this way that Theorem \ref{thmA} is essentially equivalent to the first part of \cite[Theorem 1]{DBull}, which we state now. 
\begin{thm}
\label{thmA1}
With the above notation, one has
$$\dim ER(f)_k=\defect_{n(d-2)-1-k}(\Sigma)$$
for $0 \leq k \leq n(d-2)-1$ and $\dim ER(f)_j=\tau(V)$ for $j \geq n(d-2)$.
\end{thm}
The proofs of both Theorems \ref{thmA} and \ref{thmA1} use the
Cayley-Bacharach Theorem, as discussed for instance in \cite{EGH}.
\begin{ex}
\label{exA} If we take $a=n(d-2)-d-1$, then
the hypersurface $V:f=0$ is $a$-versal if and only if the family of all hypersurfaces of degree $d$ in $\PP^n$  versally deform all the singularities of $V$, a property called {\it $T$-condition} or {\it $T$-smoothness} in \cite{GLSnb, Sh,ShT}.
This property holds if and only if
\begin{equation} 
\label{inT} 
n(d-2)-d-1< mder(f).
\end{equation} 
For instance, in the case of a plane curve, $n=2$ and the condition becomes 
$$d-5<mder(f).$$
The inequality \eqref{van2} implies that the condition \eqref{inT} holds if
$\tau(V) \leq d-1$. In fact, for $d \geq 5$, it is known that the condition \eqref{inT} holds if $\tau(V) <4(d-1)$, see \cite{GL, Sh} for the case $n=2$, and \cite{duPCTC00,ShT} for the case $n \geq 2$.
\end{ex}
One has also the following result, see \cite[Theorem 2.1]{duP1}, which we recall for the completeness of our presentation.
\begin{thm}
\label{thmB}
With the above notation,  we suppose that 
$\dim ER(f)_{a}=1$,
and $\rho=(a_0,\dots,a_n) \in ER(f)_{a}$ is a non-zero element. If there is a non-simple singular point $p \in V$ such that $\rho(p)\ne (0,\dots,0)$, then $V$ is topologically $a$-versal.
\end{thm}
\begin{ex}
\label{exB} Let $n=2$ and $V:f=x_0^d+x_1^{d-1}x_2$, with $d \geq 5$.
Then $V$ has a non-simple singularity at $p=(0:0:1)$ and $\rho=(0,x_1,-(d-1)x_2)\in ER(f)_1$ does not vanish at $p$. It follows that $V$ is topologically 1-versal.
\end{ex}

\section{Bounds on the global Tjurina number, stability and Torelli properties} \label{sec2}

One has the following result, see \cite[Theorem 5.3]{duPCTC01}.
\begin{thm}
\label{thmC}
With the above notation,  we set $r=mdr(f)$.
Then
$$(d-r-1)(d-1)^{n-1} \leq \tau(V) \leq (d-1)^n-r(d-r-1)(d-1)^{n-2}.$$
\end{thm}
For $n=2$ this result was obtained in \cite{duPCTC99}, 
and played a key role in the understanding of free curves. Indeed, when $n=2$, the reduced curve $V$ is free if and only if 
$$\tau(V)=(d-1)^2-r(d-r-1),$$
i.e. the upper bound is attained, see \cite{Dmax, E} for related results.
When $n>2$, a free hypersurface $V$ has non-isolated singularities,
and so freeness must be related to other invariants, see for instance
\cite{DmaxS}.

\begin{rk}
\label{rkC} 
The lower bound in Theorem \ref{thmC} is attained for any pair $(d,r)$.
Indeed, it is enough to find a degree $d$, reduced curve $C: f'(x_0,x_1,x_2)=0$ such that $r=mdr(f')$ and
$$\tau(C)=(d-r-1)(d-1),$$
and then take $V:f=0$, with $$f(x_0, \dots,x_n)=f'(x_0,x_1,x_2)+x_3^d+ \dots +x_n^d.$$
This formula for $f$ implies that $mdr(f)=mdr(f')$. The existence of curves $C$ as above is shown in \cite[Example 4.5]{DSt3syz}.
Similarly, the upper bound in Theorem \ref{thmC} is attained for any pair $(d,r)$ with $2r<d$, since for such pairs $(d,r)$ the existence of free plane curves $C:f'=0$ of degree $d$ and with $r=mdr(f')$ is shown in
\cite{DStExpo}. It is an interesting {\it open question} to improve the upper bound in Theorem \ref{thmC} when $2r \geq d$. The best upper bound for such pairs is (at least conjecturally) known when $n=2$, see  \cite{duPCTC99, DStMax}.
\end{rk}

 The exact sequence of coherent sheaves on $X=\PP^n$ given by
\begin{equation} 
\label{es1}
 0 \to T\langle V\rangle \to \OO_X(1)^{n+1} \to \I_{\Sigma}(d) \to 0,
\end{equation}
where the last non-zero morphism is induced by $(a_0,...,a_n) \mapsto a_0f_{0}+...a_nf_{n}$ and $\I_{\Sigma}$ is, as above, the ideal sheaf defining the singular subscheme $\Sigma$,
can be used to define the sheaf $T\langle V\rangle$ of logarithmic vector fields along $V$, see
\cite{FV1, MaVa,Se,UY}. This is a reflexive sheaf, in particular a locally free sheaf $T\langle V\rangle$ in the case $n=2$. The above exact sequence clearly yields
\begin{equation} 
\label{eqAR}
AR(f)_m=H^0(X,T\langle V\rangle(m-1)),
\end{equation} 
for any integer $m$. This equality can be used to show the reflexive sheaf $T\langle V\rangle$ is stable in many cases. This was done already in the case $n=2$ in \cite{DS14} and in the case $n=3$ in \cite{DAlgGeo}. The next result is a substantial improvement of \cite[Theorem 1.3]{DAlgGeo}.
\begin{thm}
\label{thm0}
Assume that the surface $V:f=0$ in $\PP^3$ of degree 
$d=3d'+\epsilon \geq 2$, with $d' \in \N$ and $\epsilon \in \{1,2,3\}$,
has only isolated singularities and satisfies
$$\tau(V) <(d-d'-1)(d-1)^2.$$
Then $F=T\langle V\rangle(d'-1)$ is a normalized stable rank 3 reflexive sheaf on $\PP^3$, with first Chern class
$c_1(F)=1-\epsilon \in \{0,-1,-2\}$. This reflexive sheaf is locally free if and only if $V$ is smooth.
\end{thm}

\proof
Checking the proof of \cite[Theorem 1.3]{DAlgGeo}, we see that the only point to be explained is the vanishing of $H^0(\PP^3,F)$. Using the formula \eqref{eqAR}, it follows that we have to show that $AR(f)_{d'}=0$ or, equivalently, that $r=mdr(f) >d'$. Using Theorem \ref{thmC}, we see that
$r \leq d'$ implies $\tau(V) \geq (d-d'-1)(d-1)^2$. This ends the proof of the vanishing $AR(f)_{d'}=0$.
\endproof

\begin{rk}
\label{rk0} 
The hypothesis in  \cite[Theorem 1.5]{DAlgGeo} is  
$$\tau(V) \leq 8d'+3(\epsilon-2),$$
hence the upper bound for $\tau(V)$ is, as a function of $d$, equivalent to $\frac{8}{3}d$. On the other hand, the upper bound for $\tau(V)$ in Theorem \ref{thm0} is, as a function of $d$, equivalent to
$\frac{2}{3}d^3$, hence it has a much faster growth when $d$ increases.
\end{rk}

Recall the following notion.
\begin{definition}\label{deftor}
A reduced hypersurface $V\subset \PP^n$ is called  \emph{DK-Torelli} (where DK stands for Dolgachev-Kapranov) if   the hypersurface $V$ can be reconstructed as a subset of $\PP^n$ from the sheaf
$T\langle V \rangle$.
 \end{definition}
For a discussion of this notion and various examples and results we refer to the papers \cite{DAlgGeo, DS14, DK,UY}.
The following result uses Theorem \ref{thmC} to improve \cite[Theorem 1.5]{DAlgGeo} when $n\geq 3$. More precisely we prove the following.

\begin{thm}
\label{thm1}
Let $V:f=0$ be a degree $d \geq 4$ hypersurface in $\PP^n$, having  only isolated singularities. Set $m=\lfloor \frac{d-2}{2} \rfloor$ and assume
$$\tau(V) < {m+n-1 \choose n-1}.$$
 Then one of the following holds.
\begin{enumerate}

\item $V$ is DK-Torelli;

\item $V$ is of Sebastiani-Thom type, i.e. in some linear coordinate system $(x_0,...,x_n)$ on $\PP^n$, the defining polynomial $f$ for $V$ is written as a sum $f=f'+f''$, with $f'$ (resp. $f''$) a homogeneous polynomial of degree $d$ involving only $x_0,...,x_r$ (resp. $x_{r+1},...,x_n$) for some integer $r$ satisfying $0 \leq r<n$.

\end{enumerate}

\end{thm}

\proof
We indicate only the changes to be made in the proof of \cite[Theorem 1.5]{DAlgGeo}. Let $I$ be the saturation of the ideal $J_f$ with respect to the maximal ideal $(x_0,\dots,x_n)$. The first step in the proof is to show the existence of two polynomials $h_1,h_2 \in I_m$ having no common factor. As explained in the proof of \cite[Theorem 1.5]{DAlgGeo}, to get this it is enough to assume
$$\tau(V) < {m+n-1 \choose n-1},$$
which is exactly our assumption now. The second step is to show that
$r=mdr(f) >d-2$. If we assume $r \leq d-2$, it follows from Theorem \ref{thmC} that
$$\tau(V) \geq (d-r-1)(d-1)^{n-1} \geq (d-1)^{n-1}.$$
But this is impossible, since
$${m+n-1 \choose n-1}<(d-1)^{n-1}.$$
To see this, it is enough to check that
$$\frac{m+k}{k}<d-1$$
for $k=1, \dots, n-1$ which is obvious using the definition of $m$ and the fact that $d \geq 4$.
The final step is to show that $V$ cannot have a singular point $p$ of multiplicity $d-1$. Assume such a point $p$ exists and let $g=0$ be a local equation for the singularity $(V,p)$. 
Since all the elements in $J_g+(g)$ have order at least $d-2$ and since
$$m \leq d-3,$$
the definition of $T(g)$ in \eqref{Tju} shows that the monomials in $y_j$'s of degree $m$ are linearly independent in $T(g)$. It follows that
$$\tau(V) \geq \tau(V,p)=\dim T(g) \geq {m+n-1 \choose n-1},$$
a contradiction.
\endproof

\begin{rk}
\label{rk1} 
The hypothesis in  \cite[Theorem 1.5]{DAlgGeo} is
$$\tau(V) \leq \frac{(n-1)(d-4)}{2}+1,$$
which is the same as the hypothesis above for $n=2$, but  much more restrictive for $n \geq 3$. For instance, for $n=3$ and $d=2d'$ even, the assumption in Theorem \ref{thm1} is 
$$\tau(V) < {d'+1 \choose 2},$$
while the assumption in  \cite[Theorem 1.5]{DAlgGeo} is
$$\tau(V) \leq 2d'-3.$$
\end{rk}

\end{document}